\documentclass[12pt]{amsart}
\renewcommand{\bar}{\overline}
\newcommand{\eps}{\epsilon}
\newcommand{\pa}{\partial}

\usepackage{amsmath,amsthm,amscd}
\newcommand{\frk}[1]{{\mathfrak{#1}}}

\title
[]{On the Hodge Metric of the Universal Deformation Space
of Calabi-Yau Threefolds}
\author[]{
Zhiqin Lu}
\date{July, 4th, 1997}
\address[Zhiqin Lu]
{Department of Mathematics\\
Columbia University\\
New York, NY 10027}
\email{lu@math.columbia.edu}

\newtheorem{theorem}{Theorem}[section]
\newtheorem{lemma}{Lemma}[section]

\newtheorem{prop}{Proposition}[section]

\newtheorem{definition}{Definition}[section]

\theoremstyle{remark}
\newtheorem{rem}{Remark}[section]

\begin{document}
\maketitle
\renewcommand{\baselinestretch}{1.5}

\numberwithin{equation}{section}

\section{Introductions}

A polarized Calabi-Yau manifold is a pair $(X,\omega)$ of a compact
algebraic manifold  $X$ with zero first Chern class and a K\"ahler form
$\omega\in H^2(X,Z)$. The form $\omega$ is called a  polarization. Let
${\mathcal M}$ be the universal deformation space of $(X,\omega)$.
${\mathcal M}$ is smooth by a theorem of Tian~\cite{T1}. By~\cite{Y1}, we
may assume that each $X'\in{\mathcal M}$ is a K\"ahler-Einstein manifold.
i.e.
the associated K\"ahler metric $(g'_{\alpha\bar\beta})$ is Ricci flat. The
tangent space $T_{X'}{\mathcal M}$ of ${\mathcal M}$ at $X'$ can be
identified with $H^1(X',T_{X'})_{\omega}$ where
\[
H^1(X',T_{X'})_{\omega}=\{\phi\in H^1(X',T_{X'})| \phi\lrcorner\omega=0\}
\]

The Weil-Petersson metric
$G_{PW}$ on ${\mathcal M}$ is defined by
\[
G_{WP}(\phi,\psi)=\int_{X'}
{g'}^{\alpha\bar\beta}g'_{\gamma\bar\delta}\phi^\gamma_{\bar\beta}\bar{\psi}^\delta_{\bar\alpha}
dV_{g'}
\]
where $\phi=\phi^\gamma_{\bar\beta}\frac{\partial}{\partial z^\gamma}
d\bar{z}^\beta$, 
$\psi=\psi^\delta_{\bar\alpha}\frac{\partial}{\partial z^\delta}
d\bar{z}^\alpha$ are in $H^1(X',T_{X'})_\omega$, $g'=g'_{\alpha\bar\beta}
dz^\alpha d\bar{z}^\beta$ is the K\"ahler-Einstein metric on $X$
associated with the polarization $\omega$.

In this paper, we consider the universal deformation space ${\mathcal M}$
of a simply connected Calabi-Yau threefold. Let $\omega_{WP}$ be the
K\"ahler form of the Weil-Petersson metric and set $n=\dim
H^1(X,T_X)$ for
some $X\in{\mathcal M}$. We proved

\begin{theorem}
Let $\omega_H=(n+3)\omega_{WP}+Ric(\omega_{WP})$. Then
\begin{enumerate}
\item $\omega_H$ is a K\"ahler metric on ${\mathcal M}$;
\item The holomorphic bisectional curvature of $\omega_H$ is nonpositive.
Furthermore,
Let $\alpha=((\sqrt{n}+1)^2+1)^{-1}>0$. Then the  Ricci curvature $Ric
(\omega_H)\leq
-\alpha \omega_H$ and the holomorphic sectional curvature is also less
than or
equal to $-\alpha$.
\item If $Ric(\omega_H)$ is bounded, then the Riemannian sectional
curvature of $\omega_H$ is also bounded.
\end{enumerate}
\end{theorem}

Because of the following theorem, we call $\omega_H$ the Hodge metric of 
the universal deformation space.
For the definitions, see Section 2 and Section 3. 

\begin{theorem}
Let $U$ be an open neighborhood of ${\mathcal M}$, and let $U\rightarrow D$
be the period map to  the classifying space $D$. Then
up to a constant, $\omega_H$ is the pull back of the invariant Hermitian
metric of the classifying space $D$.
\end{theorem}

\begin{rem}
In fact, we have proved  more. The theorems are also true on the
normal horizontal slices. A normal horizontal slice is a horizontal slice
such that the Weil-Petersson metric can be defined. See Section 3 for
details.
\end{rem}

The proof of the first theorem is a straightforward computation
using the Strominger's formula~\cite{S}. Using 
this
method, we can find the optimal upper bound of the Ricci curvature and the
holomorphic sectional curvature. The combination of the first and the
second
theorem is somewhat unexpected: let's explain this a little bit more in
detail. By a theorem of Griffiths, we know that the holomorphic
sectional curvature on the horizontal directions of the classifying space
is negative away from zero.
Using the same method, we know that the holomorphic bisectional curvature
are
nonpositive on certain directions. {\sl If $D$ is a homogeneous K\"ahler
manifold}, then by the Gauss theorem, we
should be able to prove
 that the holomorphic
sectional curvature and the holomorphic bisectional curvature 
of the horizontal slice are smaller than the corresponding curvatures
on the classifying space. However, $D$ is  not a homogeneous
K\"ahler manifold in general.
Nevertheless, the theorems tell us that we still have the negativity of
the curvatures.

In order to prove the second theorem, we make use of the fact that $D$ is
the dual homogeneous manifold of a K\"ahler C-space. Write $D=G/V$ where
$G$ is a noncompact semi-simple Lie group without compact factors and $V$
is its compact subgroup. Let $K$ be the maximal compact subgroup
containing $V$. We write out explicitly the projection
$G/V\rightarrow
G/K$ via local coordinate. Then the metric
$(n+3)\omega_{WP}+Ric(\omega_{WP})$ and the 
restriction of the invariant Hermitian metric of $D$ on $U$ can be
 identified.

In the last section, we gave an asymptotic estimate of the Weil-Petersson
metric to
the degeneration of Calabi-Yau threefolds.
Such an estimate was obtained by Tian~\cite{T2} in the case that the
degenerated
Calabi-Yau threefold has only ordinary double singular points. 
 C-L. Wang also got such a
result using  a completely different method.

{\bf Acknowledgment}
This paper is a refinement of a part of my ph. D  
thesis.
The author would like to thank his advisor, Professor
G. Tian for his advise and constant encouragement during my four
year's ph.D study. He also thanks Professor S. T. Yau for his constant
encouragement and many important ideas stimulating further study of this
problem. 

\section{The Classifying Space and the Horizontal Slices}

The concepts of the classifying space and the horizontal slice were
introduced by Griffiths~\cite{Gr}. We recall his definitions and notations
in this
section.

Suppose $X$ is a simply connected algebraic Calabi-Yau three-fold. The
Hodge
decomposition of the cohomology group $H=H^3(X,C)$ is
\[
H^3(X,C)=H^{3,0}\oplus H^{2,1}\oplus H^{1,2}\oplus H^{0,3}
\]
where
\[
H^{p,q}=H^q(X,\Omega^p)
\]
and $\Omega^p$ is the sheaf of the holomorphic $p$-forms.
The  quadratic form $Q$ on $X$  is defined by
\[
Q(\xi,\eta)=-\int_X\xi\wedge\eta
\]

By the Serre duality and the fact that the canonical bundle is trivial,
$\dim H^{2,1}=\dim H^{1,2}=\dim H^1(X,T_X)=n$, and 
$\dim H^{3,0}
=\dim H^{0,3}=1$. Thus $H^3(X,C)=C^{2n+2}$ is a (2n+2)-dimensional
complex vector
space.

It is easy to check that $Q$ is skew-symmetric. Furthermore, we have the
following two Hodge-Riemannian relations:

1. $Q(H^{p,q}, H^{p',q'})=0$ unless $p'=3-p$ and $q'=3-q$;

2. $(\sqrt{-1})^{p-q}Q(\psi,\bar\psi)>0$ for any nonzero element
$\psi\in H^{p,q}$.

We define the Weil operator $C:H\rightarrow H$ by
\[
C|_{H^{p,q}}=(\sqrt{-1})^{p-q}
\]
For any collection of $\{H^{p,q}\}$'s, set
\begin{align*}
& F^3=H^{3,0}\\
& F^2=H^{3,0}\oplus H^{2,1}\\
& F^1=H^{3,0}\oplus H^{2,1}\oplus H^{1,2}
\end{align*}
Then $F^1,F^2,F^3$ defines a filtration of $H$
\[
0\subset F^3\subset F^2\subset F^1\subset H
\]

Under this terminology, the Hodge-Riemannian relations can be re-written
as

3. $Q(F^3,F^1)=0, Q(F^2,F^2)=0$;

4. $Q(C\psi,\bar\psi)>0$ if $\psi\neq 0$

\vspace{0.2in}

Now we suppose that $\{h^{p,q}\}$ is a collection of integers such that
$p+q=3$ and $\sum h^{p,q}=2n+2$. 

\begin{definition}
With the notations as above, the classifying space $D$ of the Calabi-Yau
three-fold is the set of all collection of subspaces $\{ H^{p,q}\}$ of $H$
such that 
\[
H=\underset{p+q=3}{\oplus}H^{p,q}\qquad
H^{p,q}=\bar{H^{q,p}},\qquad
 \dim\,H^{p,q}=h^{p,q}
\]
and on which $Q$ satisfies the two Hodge-Riemannian relations 1,2.

Set $f^p=h^{n,0}+\cdots +h^{p,n-p}$. Then $D$ is also
the
set of all filtrations
\[
0\subset F^3\subset F^2\subset F^1\subset H, \qquad 
F^p\oplus\bar{F^{4-p}}=H
\]
with $\dim F^p=f^p$ 
on which $Q$ satisfies the bilinear relations 3,4.
\end{definition}

$D$ is a homogeneous complex manifold. 
The horizontal distribution $T_h(D)$ is defined as
\[
T_h(D)=\{X\in T(D)| XF^3\subset F^2,XF^2\subset F^1\}
\]
where $T(D)$ is the holomorphic tangent bundle which can be identified as a
subbundle of the 
(locally trivial) bundle 
$Hom(H^3(X,C),H^3(X,C))$. So $X$ naturally acts on $F^p$.

\begin{definition}
A complex integral submanifold of the horizontal distribution $T_h(D)$  is
called a horizontal slice. 
\end{definition}

Suppose $U\subset {\mathcal M}$ is a neighborhood of ${\mathcal M}$ at the
point $X$. Then there is a natural map $p: U\rightarrow D$, called the
period map, which sends a Calabi-Yau threefold to its ``Hodge Structure''.
To be precise, Let $X'\in U$. Then there is a natural identification of
$H^3(X',C)$ to $H^3(X,C)=H$. So $\{H^{p,q}(X')\}_{p+q=3}$ are the
subspaces of $H$ satisfying the Hodge-Riemannian Relations. We define
$p(X')=\{H^{p,q}(X')\}\in D$.

\section{The Weil-Petersson Metric and the Hodge Metric}
On the classifying space $D$, we can define the so called Hodge holomorphic
bundles $\underline{F}^3$, $\underline{F}^2$, $\underline{F}^1$, which are
the subbundles of the locally trivial bundle $\underline{C}^{2n+2}$. The
fiber of the bundle $\underline{C}^{2n+2}$ at $X\in{\mathcal M}$ is
$H^3(X,C)$. 
The fibers of $\underline{F}^3$, $\underline{F}^2$, $\underline{F}^1$
at $X$ are $H^{3,0}(X)$, $H^{3,0}(X)\oplus H^{2,1}(X)$,
$H^{3,0}(X)\oplus H^{2,1}(X)\oplus H^{1,2}(X)$, respectively.
Note that $\underline{F}^3$ is in fact  a line bundle. Let
$\Omega$ be a (nonzero) local holomorphic section of $\underline{F}^3$. The
curvature
form of the bundle $\underline{F}^3$ is then
$\sigma=-\frac{\sqrt{-1}}{2}\pa\bar\pa\,\log\,
Q(\Omega,
\bar\Omega)$.

Let $U$ be a horizontal slice, define
\[
\omega=\sigma|_{U}
\]

\begin{prop}\label{p31}
Let $\omega=\frac{\sqrt{-1}}2 g_{\alpha\bar\beta} dz^\alpha\wedge
d\bar{z}^\beta$ in local coordinate. Then $g_{\alpha\bar\beta}\geq 0$ is
semi-positive definite.
\end{prop}

{\bf Proof:} Let $K=-\log\,Q(\Omega,\bar\Omega)$. 
Then
\[
g_{\alpha\bar\beta}=
-\frac
{Q(\pa_\alpha\Omega+K_\alpha\Omega,
\bar{\pa_\beta\Omega+K_\beta\Omega})}
{Q(\Omega,\bar\Omega)}
\]
But $\pa_\alpha\Omega+K_\alpha\Omega\in H^{2,1}$.
The proposition follows from the second Hodge-Riemannian Relation. 
\qed

\begin{definition}
The horizontal slice is called normal if the form $\omega$ is positive
definite
at any point. In that case, $\omega=\omega_{WP}$ is called the
Weil-Petersson
metric on the normal horizontal slice.
\end{definition}

\begin{rem}
By the theorem of Tian~\cite{T1}, we know that if ${\mathcal M}$ is a
universal deformation space, then $\sigma|_{\mathcal M}$ is the
Weil-Petersson metric defined in the introduction. Thus the universal
deformation space is a normal horizontal slice.
\end{rem}

\begin{definition}
The cubic form $F=F_{ijk}$ is a (local) section of the bundle
$Sym^3(T^*{\mathcal M})\otimes (F^3)^{\otimes 2}$ defined by
\[
F_{ijk}=Q(\Omega,\pa_i\pa_j\pa_k\Omega)
\]
in local coordinates $(z^1,\cdots,z^n)$.
\end{definition}

\begin{definition}
Let $U$ be a normal horizontal slice. Suppose $\omega_D$ is the K\"ahler
form of the invariant
Hermitian metric on $D$, then we call $\omega_D|_U$ the Hodge metric on
$U$.
\end{definition}

We are going to prove

\begin{theorem}\label{31}
Suppose $\omega_{WP}$ is the
K\"ahler form of the Weil-Petersson metric. Let
\[
\omega_1=(n+3)\omega_{WP}+Ric(\omega_{WP})
\]
then $\omega_1$ is a constant multiple of the Hodge metric.
\end{theorem}

Before proving the theorem, we first prove

\begin{prop}\label{original}
There is a basis $e_1,\cdots,e_{2n+2}$ of $H$ under which $Q$
can be represented as
\[
Q=
\sqrt{-1}
\left(
\begin{array}{cc}
& 1\\
-1 &
\end{array}
\right)
\]
And if we let
\begin{align*}
&f^3=span \{e_1-\sqrt{-1}e_{n+2}\}\\
&f^2=span\{e_1-\sqrt{-1}e_{n+2},e_2+\sqrt{-1}e_{n+3},
\cdots, e_{n+1}+\sqrt{-1}e_{2n+2}\}
\end{align*}
and $f^1$ is the hyperplane perpendicular to $f^3$ with respect
to $Q$, then
\[
\{0\subset f^3\subset f^2\subset f^1\subset H\}\in D
\]
\end{prop}\qed

 The point $\{f^3,f^2,f^1\}\in D$ is called the original point of $D$.
Sometimes we write it as $eV$ if $D=G/V$.

By the curvature formula of Strominger~\cite{S}, 
the Ricci curvature of the Weil-Petersson metric
\[
R_{i\bar j}=-(n+1)g_{i\bar j}+e^{2K}F_{ipq}\bar{F_{jmn}}g^{p\bar m}
g^{q\bar n}
\]
where we set $\omega_{WP}=\frac{\sqrt{-1}}{2}g_{i\bar j}dz^i\wedge
d\bar{z}^j$ in
the
local coordinates
 $(z^1,\cdots, z^n)$
and $K$ is the local function: $K=-\log Q(\Omega,\bar\Omega)$.

 Let
$\omega_1=\frac{\sqrt{-1}}{2}
h_{i\bar j}dz^i\wedge d\bar{z}^j$. Then 
\[
h_{i\bar j}=2g_{i\bar j}+e^{2K}F_{ipq}\bar{F_{jmn}}
g^{p\bar m}g^{q\bar n}
\]
Clearly $(h_{i\bar j})>0$.

Suppose
\[
\pi': D\rightarrow CP^{2n+1}
\]
is the projection of $D$ to $CP^{2n+1}$ by sending
$(F^3,F^2,F^1)$ to $F^3$. 
Let $U$ be a normal horizontal slice. 
Let $\Omega$ be a (nonzero) local section of $\underline{F^3}$. Then
$\pa_i\Omega+K_i\Omega$ is not zero because $\omega_{WP}$ is positive. Thus
\[
\pi': U\rightarrow D\rightarrow CP^{2n+1}
\]
is an immersion.

Now we consider the result of  Bryant and Griffiths~\cite{BG}. Their
results can
be
briefly
written as follows:

We assume that $eV\in U$. i.e. the normal horizontal  slice
passes the original point of $D$, where the original point is defined
as $\{f^3,f^2,f^1\}\in D$ in the Proposition~\ref{original}.
Then according to Bryant and Griffiths, there is a holomorphic
function $u$ defined on a neighborhood of the original point of 
${\mathbb C}^n$ such
if $(z^1,\cdots, z^n)$ is the local holomorphic coordinate of $U$ at $eV$,
the original point, then
\begin{equation}\label{o}
\Omega=(1,\frac{1}{\sqrt{2}}z^1,\cdots,\frac{1}{\sqrt{2}}z^n,
u-\sum_i\frac 12z^iu_i,\frac{1}{\sqrt{2}}u_1,\cdots,\frac{1}{\sqrt{2}}u_n)
\end{equation}
with  $F^1=Span \{ (\Omega)\}$,  $F^2=span \{\nabla \Omega\}$, $F^1\perp
F^3$ via $Q$. In particular, $u(0)=-\sqrt{-1}$, $|\nabla u(0)|=0$,
$\nabla^2
u(0)=\sqrt{-1}I$, where $I$ is the unit matrix. 

In order to prove Theorem~\ref{31}, we need only to prove it at the
original 
point, because any point of the homogeneous space $D$ can be
taken as the original point.

Under these notations, at $eV$,  
$D_i\Omega=\pa_i\Omega+K_i\Omega=\frac{1}{\sqrt{2}}(e_{i+1}+\sqrt{-1}e_{n+i+2})$.
Thus

\[
\sqrt{-1}Q(\Omega,\bar\Omega)=-2
\]
and
\[
g_{i\bar j}=-\frac{Q(D_i\Omega,\bar{D_j\Omega})}
{Q(\Omega,\bar{\Omega})}=\frac 12\delta_{ij}
\]
Furthermore, the cubic form $F_{ijk}$ at $eV$ is 
\[
F_{ijk}=-\frac 12\frac{\pa^3 u}{\pa z^i\pa z^j\pa z^k}(0)
=-\frac 12 u_{ijk}(0)
\]
Thus 
\[
h_{i\bar j}=2g_{i\bar j}+e^{2K}F_{imn}\bar{F_{jpq}}g^{m\bar p}g^{n\bar q}
=\delta_{ij}+\frac 14 u_{imn}(0)\bar{u_{jmn}}(0)
\]

Now we are going to prove that $(h_{i\bar j})$ is a constant multiple of
the Hodge metric.
Consider the projection
\[
\pi: D=G/V\rightarrow G/K
\]
where $K$ is the maximal
connected compact subgroup of $G$ containing $V$.
We have

\begin{lemma}
Let $U$ be a horizontal slice, then
$\pi$ is an isometry between the Riemannian submanifold $U$
of $D$
 and the Riemannian submanifold $\pi(U)$ of $G/K$.
\end{lemma}

{\bf Proof:} Note that $U$ is a horizontal slice of $D$. The
lemma follows from the definition of the
invariant Hermitian metric on both
manifold.

\qed

From the above lemma, we know that in order to compute $\omega_D|_U$, we
need only computed the metric of $U$ as a submanifold of $G/K$, even the
map $\pi$ is not holomorphic(Recall that $D$ is not homogeneous K\"ahler,
so the map will not be holomorphic in general).
In order to do this, 
we write out  the projection 
\[
\pi : G/V\rightarrow G/K
\]
explicitly now.

\vspace{0.15in}

It is easy to prove from linear algebra that 
 the projection $\pi$ send
\[
0\subset F^3\subset F^2\subset F^1\subset H
\]
to
\[
F^3\oplus H^{1,2}
\]

We have known that $G/K=Sp(n+1, R)/U(n+1)$ is the Hermitian symmetric
space. $G/K$ can be realized as the set of $(n+1)$ planes $P$ in the
$C^{2n+2}$ space such that $-\sqrt{-1}Q(P,\bar{P})>0$. Thus $G/K$ can be
represented as
the set of all the symmetric $(n+1)\times(n+1)$ matrix $Z$
satisfying
$Im\,Z>0$
where $Im\,Z>0$ means $Im\,Z$ is a positive definite Hermitian matrix.

We  write  the entries of the matrix $Z$ as  functions of $D$.

Suppose now 
near the original point, $F^3$ is spanned by
\[
(1,z^t,a,\alpha^t)
\]
where $z,\alpha\in {\mathbb C}^n, a\in {\mathbb C}$. And suppose $F^2$ 
is spanned by the row vectors of the matrix
\[
\left(
\begin{array}{llll}
1  &  z^t   & a   &  \alpha^t\\
0   & 1   &  \beta   &   A
\end{array}
\right)
\]
for $\beta\in {\mathbb C}^n, A\in {\frk g}{\frk l}(n,{\mathbb C})$. 
Then by the first Hodge-Riemannian relation
\[
Q(F^2,F^2)=0
\]
we know that
\[
\beta=\alpha-Az,\qquad A^t=A
\]
So locally, we can represented $F^2$ by the matrix
\begin{equation}
\label{eq25}
\left(
\begin{array}{cccc}
1   &  z^t  & a  &  \alpha^t\\ 
0  &  1   &  \alpha-Az  &A
\end{array}
\right)
\qquad A^t=A
\end{equation}

Let 
$\Omega=(1,z^t,a,\alpha^t)$, and let
$\binom{\Omega}{\Theta}$ be a local section of $F^2$ with
$\Theta=(0,1,\alpha-Az,A)$. Set  
\begin{align*}
&m=Q(\Omega,\bar\Omega)=-a+\bar{a}-\alpha^t\bar{z}+\bar{\alpha}^tz\\
&\xi=Q(\Omega,\bar\Omega)=-\alpha+\bar{\alpha}-\bar{A}(\bar{z}-z)
\end{align*}
where $m\in{\mathbb C}$, $\xi\in{\mathbb C}^n$. 
It is easily checked that
\[
Q(\Omega,\bar{\Theta}-\frac{\xi}{m}\bar\Omega)=0
\]

So $\bar\Theta-\frac{\xi}{m}\bar{\Omega}\in F^1$ and since $\Omega$ and 
$\Theta$ are in $H^{2,1}$, $\bar\Theta-\frac{\xi}{m}\bar{\Omega}\in
H^{1,2}$.

The projection $\pi$ can be locally written as
\begin{align*}
&\left(
\begin{array}{cccc}
1 & z^t  &  a  & \alpha^t\\
0  &  1  &  \alpha-Az  &  A  
\end{array}
\right)
\longrightarrow
\left(
\begin{array}{c}
\Omega\\
\bar\Theta-\frac{\xi}{m}\bar\Omega
\end{array}
\right)=\\
&\left(
\begin{array}{cccc}
1   &  z^t   &  a   &  \alpha^t\\
-\frac{\xi}{m}  &  1-\frac{\xi\bar{z}^t}{m}   & 
\overline{\alpha-Az}-\frac{\bar{a}\xi}{m}  &
\bar{A}-\frac{\xi}{m}\bar\alpha^t
\end{array}
\right)
\end{align*}
where the right hand side of the above represents an $(n+1)$-plane in $H$.

The symmetric $(n+1)\times (n+1)$ matrix $Z$ can be obtained as follows:
let
\[
\mu=\frac{1}{m-(\bar{z}^t-z^t)\xi}
\]

Then as a matrix
\[
1+\mu\xi(\bar{z}^t-z^t)=(1-\frac{\xi}{m}(\bar{z}^t-z^t))^{-1}
\]
We have

\[
\left(
\begin{array}{cc}
1  &  z^t\\
-\frac{\xi}{m}  &  1-\frac{\xi\bar{z}^t}{m}
\end{array}
\right)
^{-1}
=
\left(
\begin{array}{cc}
1-z^tB\frac{\xi}{m}   &   -z^tB\\
B\frac{\xi}{m}     &  B
\end{array}
\right)
\]
where $B=1+\mu\xi(\bar{z}^t-z^t)$. Let
\[
\left(
\begin{array}{cc}
1-z^tB\frac{\xi}{m}   &   -z^tB\\
B\frac{\xi}{m}     &  B
\end{array}
\right)
\left( 
\begin{array}{cc}
a   &  \alpha^t\\
\overline{\alpha-Az}-\frac{\bar{a}\xi}{m}  &
\bar{A}-\frac{\xi}{m}\bar\alpha^t
\end{array}
\right)
=
\left( 
\begin{array}{cc}
D_1   &   D_2^t\\
D_3   &  D_4
\end{array}
\right)
\]
for $D_1\in {\mathbb C}, D_2,D_3\in {\mathbb C}^n, D_4\in 
{\frk g}{\frk l}(n,{\mathbb C})$.
Then it can be computed
\begin{equation}\label{db}
\left\{
\begin{array}{l} 
D_1=a-z^t(\alpha-\bar{A}z)+\mu (z^t\xi)^2\\
D_2=D_3=(a-\bar{a})\mu\xi+B\overline{(\alpha -Az)}\\
D_4=\bar{A}+\mu\xi\xi^t
\end{array}
\right.
\end{equation}

Then the matrix $Z$ is obtained:

\begin{prop}
Under the notation as above, the map
\[
\pi: G/V\rightarrow G/K
\]
under the local coordinate described as above is
\[
\left( 
\begin{array}{cccc}
1   &   z^t  &  a  &  \alpha^t\\
    &  1    & \alpha-Az   & A
\end{array}
\right)
\longrightarrow
\left(
\begin{array}{ll}
D_1   &   D^t_2\\
D_3    &   D_4
\end{array} 
\right)
=Z
\]
where the $D_i$'s are defined as in equation~\eqref{db}.
\end{prop}

The Hermitian metric on $G/K$ is $-\sqrt{-1}\pa\bar\pa\log\det\,Im\, Z$.
In particular, at the original point, it is
$\sum_{ij}d\,Z^{ij}\wedge d\,\bar{Z^{ij}}$, where we set $Z^{ij}=Z^{ji}$
if
$i>j$. 
By Equation~\ref{o} we see that 
in order to get the map $U\rightarrow G/K$,
$z,a,\alpha, A$ in Equation~\ref{db}
should be replaced by $\frac{1}{\sqrt{2}}z, u-\sum\frac 12 z^iu_i,
\frac{1}{\sqrt{2}}\nabla u, \frac{1}{\sqrt{2}}u_{ij}$, respectively.
Thus we have
\begin{equation}
\label{eq3}
\left\{
\begin{array}{l}
\frac{\partial D_1}{\partial t_k}=0,\frac{\partial D_1}{\partial 
\bar{t}_k}=0\\
\frac{\partial (D_3)_r}{\partial t_k}=-\sqrt{2}i\delta_{rk},
\frac{\partial(D_3)_r}{\partial \bar{t}_k}=0\\
\frac{\partial D_4}{\partial t_k}=0, 
\frac{\partial (D_4)_{rs}}{\partial\bar{t}_k}=\bar{u}_{rsk}
\end{array}
\right.
\end{equation}
By a straightforward
computation, we know the restriction of the metric on $G/K$ on $U$ at the
original point is a constant multiple of
\[
h_{i\bar j}=\delta_{ij}+\frac 14 u_{imn}(0)\bar{u_{jmn}}(0)
\]
Thus completes the proof. 

\section{The Curvature Computation}\label{se:5}

In this section we give
an optimal estimate of the upper bound of 
the holomorphic sectional curvature,
bisectional curvature and the Ricci curvature
of  a normal horizontal slice. 

Let $U$ be a normal horizontal slice. 
Suppose $(g_{i\bar j})$ is the Weil-Petersson metric, $(F_{ijk})$ is the
cubic form, and $K=-\log Q(\Omega,\bar\Omega)$. 
The Hodge  metric $(h_{i\bar j})$ is:
\[
h_{i\bar j}=2g_{i\bar j}+\sum_{rspq}e^{2K}F_{irs}\bar{F_{jpq}}
g^{r\bar p}g^{s\bar q}
\]
As we have proved, $(h_{i\bar j})$ is a K\"ahler metric. So 
 the curvature tensor
$\tilde{R}_{i\bar jk\bar l}$ of $(h_{i\bar j})$ is

\[
\tilde{R}_{i\bar jk\bar l}=
\frac{\pa^2 h_{i\bar j}}{\pa z^k\pa\bar{z}^l}
-h^{n\bar m}\frac{\pa h_{i\bar m}}{\pa z^k}
\frac{\pa h_{n\bar j}}{\pa\bar{z}^l}
\]
Now we suppose that at point $p$, the local coordinate 
for the Weil-Petersson metric is normal. i.e., at
point $p$, $g_{i\bar j}=\delta_{ij}$ and $dg_{i\bar j}=0$.
Furthermore, assume $K(p)=0$. 
 The curvature
tensor $R_{i\bar jk\bar l}$ of $(g_{i\bar j})$
then is 
\[
R_{i\bar jk\bar l}=
\frac{\pa^2g_{i\bar j}}{\pa z^k\pa\bar{z}^l}
\]
also we have
\begin{equation}\label{oo}
\frac{\pa h_{i\bar m}}{\pa z^k}
=\sum_{rs}F_{irs,k}\bar{F_{mrs}}
\end{equation}
where
\[
F_{irs,k}=\pa_k F_{irs}+2K_k F_{irs}
\]
is the covariant derivative 
of the cubic form with respect to the Weil-Petersson metric.
By using the Strominger formula at $p$,
\[
R_{i\bar jk\bar l}=\delta_{ij}\delta_{kl}
+\delta_{il}\delta_{kj}-F_{ikm}\bar{F_{jlm}}
\]
We get
\begin{align}\label{y}
\begin{split}
& \frac{\pa^2h_{i\bar j}}{\pa z^k\pa\bar{z}^l}=
2R_{i\bar jk\bar l}-2\sum_{sqr}R_{q\bar sk\bar
l}F_{irs}\bar{F_{jrq}}\\
&+2\delta_{kl}\sum_{rs}F_{irs}\bar{F_{jrs}}
+\sum_{rs}F_{irs,k}\bar{F_{jrs,l}}
\end{split}
\end{align}

Combining Equation~\ref{oo} and Equation~\ref{y},we have

\begin{prop}\label{prop:p1}
If $K=0$ at the point $p$, 
\begin{align}\label{eq:eq21}
\begin{split}
&\tilde{R}_{i\bar jk\bar l}=
2R_{i\bar jk\bar l}+2\delta_{kl}\sum_{rs}F_{irs}\bar{F_{jrs}}
-2\sum_{sqr}R_{q\bar sk\bar l}F_{irs}\bar{F_{jrq}}\\
&+\sum_{rs}F_{irs,k}\bar{F_{jrs,l}}
-\sum_{mn}(\sum_{rs}F_{irs,k}\bar{F_{mrs}})
\bar{(\sum_{rs}F_{jrs,l}\bar{F_{nrs}})}h^{n\bar m}
\end{split}
\end{align}
\end{prop}

Based on the above proposition, we get
\begin{theorem}
Let $c(n)=((\sqrt{n}+1)^2+1)$, then
\begin{align*}
&Ric(\omega_H)\leq -\frac{1}{c(n)}\omega_H\\
&R\leq -\frac{1}{c(n)}
\end{align*}
where $R$ is the superium of the holomorphic sectional curvature.
The constant here is optimal. 
Furthermore, the bisectional curvature is nonpositive.
\end{theorem}

{\bf Proof:}
We  consider the point $p$ and the normal coordinate at $p$ with
respect to the Weil-Petersson metric.
Fixing $i$, let
\[
A_m=\sum_{rsk}F_{irs,k}a_k\bar{F_{mrs}}
\]
for a vector $a=(a_1,\cdots,a_n)$.
Then it is easy to see that
\begin{align}\label{eq:eq20}
\begin{split}
&\sum_{pq}|\sum_kF_{ipq,k}a_k|^2
-\sum_{mn}(\sum_{rsk}F_{irs,k}a_k\bar{F_{mrs}})
\bar{(\sum_{rsk}F_{irs,k}a_k\bar{F_{nrs}})}h^{n\bar m}\\
&=\sum_{pq}
|\sum_kF_{ipq,k}a_k-\sum_{mn}h^{n\bar m}A_mF_{npq}|^2
+2\sum_a|\sum_mh^{a\bar m}A_m|^2
\end{split}
\end{align}
where we use the fact that $h_{i\bar j}=2\delta_{ij}
+F_{imn}\bar{F_{jmn}}$ at $p$.

Define a generic vector $a^k=\delta_{ik},k=1,\cdots,n$.
Using Equation~\ref{eq:eq20}, we have
\[
\sum_{pq}|F_{ipq,i}|^2-\sum_{mn}h^{n\bar m}
(\sum_{rs}F_{irs,i}\bar{F_{mrs}})
\bar{\sum_{rs}F_{irs,i}\bar{F_{nrs}})}\geq 0
\]

Now using Proposition~\ref{prop:p1}, we get
\begin{align*}
&\tilde{R}_{\alpha\bar\beta\gamma\bar\delta}a^\alpha\bar{a^\beta}a^\gamma\bar{a^\delta}=
\tilde{R}_{i\bar ii\bar i}\geq 2R_{i\bar ii\bar i}
+2\sum_{rs}|F_{irs}|^2-2\sum_{sqr}R_{q\bar si\bar i}F_{irs}
\bar{F_{irq}}\\
&\geq 4-4\sum_r|F_{iir}|^2+2\sum_{rp}
|\sum_qF_{qip}\bar{F_{irq}}|^2
\end{align*}

Let 
\[
x=\sum_r|F_{iir}|^2
\]
Then we have
\begin{align*}
&\sum_{rp}|\sum_qF_{qip}\bar{F_{irq}}|^2\geq
|\sum_q|F_{iiq}|^2|^2=x^2\\
&\sum_{rp}|\sum_qF_{qip}\bar{F_{irq}}|^2\geq
\sum_r|\sum_q|F_{qir}|^2|^2\geq\frac{1}{n}(h_{i\bar i}-2)^2
\end{align*}

So for $a,b>0, a+b=1$, we have
\begin{align*}
&\frac 12\tilde{R}_{i\bar ii\bar i}\geq
2-2x+ax^2+\frac bn(h_{i\bar i}-2)^2\\
&\geq 2-\frac 1a+\frac bn(h_{i\bar i}-2)^2
\end{align*}

On the other hand, we have
\[
h_{\alpha\bar\beta}a^\alpha\bar{a^\beta}=h_{i\bar i}
\]

Let $a=\frac{2+\sqrt{n}}{2+2\sqrt{n}}$ and $b=1-a$, we have

\[
\tilde{R}_{i\bar ii\bar i}\geq \frac {1}{(\sqrt{n}+1)^2+1}
h_{i\bar i}^2=\frac{1}{c(n)}|h_{i\bar j}a^i\bar{a^j}|^2
\]

It is a straightforward computation that the constant here is optimal. 

Thus we proved $\tilde{R}(a,\bar a,a,\bar a)\geq ||a||^2$. Since $a$ can
be any vector by making a linear transformation of the normal coordinate.
We have already proved the assertion of the theorem about the holomorphic
sectional curvature.

Now we turn to the bisectional curvature. For any
$(a^1,\cdots, a^n)$, using the same inequalities before, we have
\begin{align*}
&\tilde{R}_{i\bar ik\bar l}a^k\bar{a^l}
\geq 2\sum_k|a^k|^2+2|a^i|^2-4\sum_r|\sum_kF_{irk}a^k|^2\\
&+2\sum_{mr}|\sum_{qk}F_{qkm}\bar{F_{irq}}a^k|^2
\geq 2|a^i|^2+2\sum_r
|\sum_{qk}F_{iqk}\bar{F_{iqr}}a^k-a^r|^2\geq 0
\end{align*}

This proves the nonpositivity of the bisectional curvature.

Finally  we consider the Ricci curvature.  Suppose that $\xi$ is a unit
vector.
Then by the definition of the Ricci curvature and above results, we have
\[
-Ric(\xi,\bar\xi)\geq \tilde{R}(\xi,\bar\xi,\xi,\bar\xi)
\]
This completes the proof of the theorem.

\section{The Boundness of the Sectional Curvature}

In this section, we prove that the boundness of the Ricci curvature
implies the boundness of the Riemannian sectional curvature. 

\begin{theorem}\label{bdd}
Suppose $U$ is a normal horizontal slice. 
Suppose $p\in U$ is a fixed point such that the Ricci curvature has a
lower bound $C_p$ at $p$. That is
\[
Ric(\omega_H)_p\geq -C_p(\omega_H)_p
\]
Then the Riemannian sectional curvature has a bound
\[
|\tilde R(X,Y,X,Y)|\leq (3+C_p) ||X||^2||Y||^2
\]
where $X,Y\in T_pU$ and $X\perp Y$.
\end{theorem}

We begin by restating Proposition~\ref{prop:p1} in the
section~\ref{se:5}.

\begin{prop}\label{prop:p22}
Suppose we have the notations as in the
proposition~\ref{prop:p1}, then we have
\[
\tilde{R}_{i\bar{j}k\bar{l}}=A_{i\bar{j}k\bar{l}}
+B_{i\bar{j}k\bar{l}}
\]
where
\begin{align*}
&
A_{i\bar{j}k\bar{l}}=2\delta_{ij}\delta_{kl}
+2\delta_{il}\delta_{kj}-4\sum_sF_{iks}\bar{F_{jls}}
+2\sum_{mnpq}F_{qkm}\bar{F_{plm}}F_{inp}\bar{F_{jnq}}\\
&
B_{i\bar{j}k\bar{l}}=\sum_{rs}
(F_{irs,k}-\sum_{mn}A_{ik\bar m}F_{nrs}h^{n\bar m})
\bar{(F_{jrs,l}-\sum_{mn}A_{jl\bar m}F_{nrs}h^{n\bar m})}\\
&
+2\sum_{mm_1n}A_{ik\bar m}h^{n\bar m}\bar{\bar{A_{jl\bar m_1}}h^{m_1\bar
n}}
\end{align*}
here we define
\[
A_{ik\bar m}=\sum_{rs}F_{irs,k}\bar{F_{mrs}}
\]
\end{prop}

{\bf Proof:} A straightforward computation. \qed

\begin{lemma}
Suppose that $\xi,\eta\in T_p^{(1,0)}{U}$ and define
$||\xi||^2=h_{i\bar j}\xi^i\bar{\xi^j}$, then
\[
|\tilde{R}_{i\bar jk\bar l}\xi^i\xi^k\bar{\eta^j}\bar{\eta^l}|
\leq (6+C_p)||\xi||^2||\eta||^2
\]
\end{lemma}

{\bf Proof:} Note that the holomorphic bisectional curvature of ${U}$ 
is
nonpositive. We know that the holomorphic sectional curvature is
bounded by $C_p$, i.e.
\[
|\tilde{R}_{i\bar jk\bar l}\xi^i\xi^k\bar{\xi^j}\bar{\xi^l}|
\leq C_p||\xi||^4
\]

We have
\begin{align*}
&\sum_{pq}|\sum_{nij}F_{inp}\bar{F_{jnq}}\xi^i\bar{\eta^j}|^2\\
&=\sum_{pq}|\sum_n(\sum_iF_{inp}\xi^i)\bar{
(\sum_jF_{jnq}\eta^j})|^2\\
&\leq\sum_{pq}(\sum_n|\sum_iF_{inp}\xi^i|^2
\sum_n|\sum_jF_{jnq}\eta^j|^2)\\
&=\sum_{pn}|\sum_iF_{inp}\xi^i|^2
\sum_{qn}|\sum_jF_{jnq}\eta^j|^2\\
&\leq ||\xi||^2||\eta||^2
\end{align*}

here we use the fact that
\[
h_{i\bar j}=2\delta_{ij}+\sum_{rs}F_{irs}\bar{F_{jrs}}
\]

and

\[
\sum_{pn}|\sum_iF_{inp}\xi^i|^2=
\sum_{i,j}\sum_{pn}F_{inp}\xi^i
\bar{F_{jnp}\xi^j}\leq||\xi||^2
\]

Thus
\begin{align*}
&\sum_{ijkl}\sum_{mnpq}F_{qkm}\bar{F_{plm}}F_{inp}
\bar{F_{jnq}}\xi^i\xi^k\bar{\eta^j}\bar{\eta^l}\\
&=\sum_{pq}(\sum_{nij}F_{inp}\bar{F_{jnq}}\xi^i\bar{\eta^j})
(\sum_{mkl}F_{qkm}\bar{F_{plm}}\xi^k\bar{\eta^l})\\
&\leq\sqrt{
\sum_{pq}|\sum_{nij}F_{inp}\bar{F_{jnq}}\xi^i\bar{\eta^j}|^2
\sum_{pq}|\sum_{mkl}F_{qkm}\bar{F_{plm}}\xi^k\bar{\eta^l}|^2}\\
&\leq ||\xi||^2||\eta||^2
\end{align*}

We also have

\[
\sum_s|\sum_{ik}F_{iks}\xi^i\xi^k|^2\leq
\sum_s(\sum_k|\sum_iF_{iks}\xi^i|^2\sum_k|\xi^k|^2)
\leq ||\xi||^4
\]

Thus by proposition~\ref{prop:p22}, we have
\[
|\sum_{ijkl}A_{i\bar jk\bar l}\xi^i\xi^k\bar{\eta^j}\bar{\eta^l}|
\leq 6||\xi||^2||\eta||^2
\]

We also have
\[
\sum_{ijkl}B_{i\bar jk\bar l}\xi^i\xi^k\bar{\eta^j}\bar{\eta^l}
\leq
\sqrt{
\sum_{ijkl}B_{i\bar jk\bar l}\xi^i\xi^k\bar{\xi^j}\bar{\xi^l}
\sum_{ijkl}B_{i\bar jk\bar l}\eta^i\eta^k\bar{\eta^j}\bar{\eta^l}}
\]

Combining the above two inequalities we proved the lemma.\qed

{\bf Proof of Theorem~\ref{bdd}:}

Let 
\begin{align*}
&\xi=X-\sqrt{-1}JX\\
&\eta=Y-\sqrt{-1}JY
\end{align*}

Then
\[
\tilde R(X,Y,X,Y)=\frac 18
(Re\,\tilde R(\xi,\bar\eta,\xi,\bar\eta)
-\tilde R(\xi,\bar\xi,\eta,\bar\eta))
\]

The bisectional curvature is bounded by $C_p$
\[
|\tilde R(\xi,\bar\xi,\eta,\bar\eta)|\leq C_p||\xi||^2||\eta||^2
\]

Thus
\[
|\tilde{R}(X,Y,X,Y)|\leq\frac 14(3+C_p)
||\xi||^2||\eta||^2=(3+C_p)||X||^2
||Y||^2
\]

\section{An Asymptotic Estimate}
In this section, we make use of the results in the previous sections
to prove an asymptotic estimate of the Weil-Petersson metric of
the degeneration of Calabi-Yau threefolds. The motivation
for the estimate is from a result of G. Tian~\cite{T2}. Although the 
argument can be generalized to study the Weil-Petersson
metric of a normal horizontal slice near infinity, we restrict
ourselves to the degeneration of Calabi-Yau threefolds.

We say $\pi:{\frk X}\rightarrow\Delta$ is a degeneration
of Calabi-Yau threefolds, if ${\frk X}$, $\Delta$ are complex 
manifolds and $\pi$ is holomorphic, and $\Delta$ is the unit disk
in ${\mathbb C}$. $\forall t\in\Delta$, $t\neq 0$,
$\pi^{-1}(t)$ is a smooth Calabi-Yau threefold while $\pi^{-1}(0)$
is a divisor of normal crossing. We also denote $\Delta^*$ to be
the punctured unit disk.

\begin{theorem}
Suppose ${\frk X}\rightarrow\Delta$ is a degeneration of Calabi-Yau
threefolds. Suppose $\omega$ be the
Weil-Petersson metric on $\Delta^*$. Then if
\[
\underset{r\rightarrow 0}{\lim}\frac{\log\,\omega}{\log \frac{1}{r}}=0
\]
Then
\[
\omega\leq \sqrt{-1}C_1(\log \frac{1}{r})^{4c(n)}dz\wedge d\bar z
\]
where $c(n)=((\sqrt{n}+1)^2+1)$, $C_1$ is a constant and $z$ is the
coordinate of $\Delta$.
\end{theorem}

\begin{rem}
$\omega$ is a K\"ahler metric on $\Delta^*$. So there is a function
$\lambda(z)>0$ on $\Delta^*$ such that
\[
\omega=\sqrt{-1}\lambda(z)dz\wedge d\bar z
\]
The assumption is understood as
\[
\underset{r\rightarrow 0}{\lim}\frac{\log \lambda}{\log\frac 1r}=0
\]
and the  
conclusion of the theorem is understood as 
\[
\lambda(z)\leq C_1(\log\frac{1}{r})^{4c(n)}
\]
\end{rem}

{\bf Proof:} By a theorem of Tian, there are no obstructions towards the
deformation of a Calabi-Yau three-fold. Suppose $M\in\pi^{-1}(\Delta^*)$
is a fiber. Let 
$n=\dim H^1(M,\Theta)$ and $\pi(M)=p$. Let
${\mathcal M}$ be the universal deformation space at $p$. Then there is a
neighborhood $U$ of $\Delta^*$ at $p$ such that $U\subset{\mathcal M}$
Suppose $z_1=z$, 
and suppose $U$ is defined
by $z_2=\cdots =z_n=0$ near the
point $p$. 
Using the notations as in the previous sections, 
by the Strominger's formula, we have
\begin{equation}\label{eq:eq56}
R_{1\bar 11\bar 1}=2g_{1\bar 1}^2-\frac{1}{(\Omega,\bar\Omega)^2}
F_{11\xi}\bar{F_{11\eta}}g^{\xi\bar \eta}
\end{equation}
Here $(\Omega,\bar\Omega)=\sqrt{-1}Q(\Omega,\bar\Omega)$

\begin{lemma}
Let $\lambda=g_{1\bar 1}$, then
\[
\lambda^{-1}\frac{1}{(\Omega,\bar\Omega)^2}
F_{11\xi}\bar{F_{11\eta}}g^{\xi\bar \eta}\leq
\frac{1}{(\Omega,\bar\Omega)^2}
F_{1\alpha\xi}\bar{F_{1\beta\eta}}g^{\alpha\bar\beta}
g^{\xi\bar\eta}
\]
\end{lemma}

{\bf Proof:}   If
$g_{1\bar\beta}=0$
for $\beta\neq 1$, then the inequality is trivially true. Thus we would
like
to
choose a coordinate such that $g_{1\bar\beta}=0,\beta\neq 1$.

Let $A$ be an $(n-1)\times (n-1)$ matrix. Let
\[
w_i=\sum_{j=2}^nA_{ij}z_j
\]
for $i=2,\cdots,n$. If $A$ is a nonsingular matrix, then
$(z_1,w_2,\cdots,w_n)$ will be  local holomorphic coordinate 
of ${\mathcal M}$ at $p$
and
$\Delta^*$ is again be defined by $w_2=\cdots=w_n=0$. Now we choose an
$A$ such that
\[
(\frac{\pa}{\pa z_1},
\frac{\pa}{\pa\bar{w}^k})=0
\]
for $k=2,\cdots, n$. Suppose $\tilde{g}_{\alpha\bar\beta}$ is the matrix
under the coordinate $(z_1,w_2\cdots,w_n)$, $(\tilde{g}_{1\alpha})=0$ for
$\alpha\neq 1$.

\qed

Using the lemma, from Equation~\ref{eq:eq56}, we have
\[
R_{1\bar 11\bar 1}\geq
2\lambda^2-\lambda(h_{1\bar 1}-2g_{1\bar 1})
\geq -\lambda h_{1\bar 1}
\]
Suppose $\Delta=\frac{\pa^2}{\pa z_1{\bar\pa z_1}}$, then the Gauss
curvature of $\Delta^*$ with respect to $\lambda$ is
\[
K=-\frac{4}{\lambda}\Delta\log\lambda
\]
On the other hand
\[
\Delta\log\log\frac 1r=-\frac{1}{4r^2(\log\frac 1r)^2}
\]

So we have
\[
\Delta\log\lambda=-\frac 14\lambda K\geq\frac{1}{\lambda}
R_{1\bar 11\bar 1}\geq -h_{1\bar 1}
\]
where we use the Gauss formula $4R_{1\bar 11\bar 1}\leq -K\lambda^2$.
The holomorphic sectional curvature of $(h_{1\bar 1})$
is less than $-\frac{1}{c(n)}$. Thus by the Schwartz lemma
\[
h_{1\bar 1}\leq\frac{c(n)}{(r\log\frac 1r)^2}
\]
Thus
\[
\Delta\log\frac{\lambda}{(\log\frac 1r)^{4c(n)}}
\geq0
\]

The rest of the proof is quite elementary: let
\[
f=\log\frac{\lambda}{(\log\frac 1r)^{4c(n)}}
\]
Then
\[
\underset{r\rightarrow 0}{\lim}\frac{f}{\log\frac 1r}=0
\]
So for any $\eps$, there is a $\delta$ such that $r<\delta$ implies
$-f+\eps\log\frac 1r$ large enough. Now
\[
\Delta(-f+\eps\log\frac 1r)\leq 0
\]
So the minimum point must be obtained at $|r|=\frac 12$. 
So
\[
f-\eps\log\frac 1r\leq C_1+\eps\log 2\leq 2C_1
\]
for any $\eps$ small. thus letting $\eps\rightarrow 0$, we have
$f\leq 2C_1$, which completes the proof.

\begin{rem}
In Hayakawa~\cite{Ha}, the author claimed a relation
between the degeneration of the Calabi-Yau manifolds
and the noncompleteness of the Weil-Petersson metric.
But her proof was incomplete. C-L. Wang~\cite{Wong} gave a  
proof of this and studied the Weil-Petersson metric in great detail.
In particular, he  proved an asymptotic estimate for the degeneration
of Calabi-Yau manifold which is slightly sharper then our estimate
independently using a different method.
\end{rem}

\bibliographystyle{plain}
\bibliography{bib}

\begin{thebibliography}{1}

\bibitem{BG}
Robert Bryant and Phillip Griffiths.
\newblock {Some Observations on the Infinitesimal Period Relations for Regular
  Threefolds with Trivial Canonical Bundle}.
\newblock In Michael Artin and John Tate, editors, {\em Arithmetic and
  Geometry}, pages 77--85. Boston, Birkha\"user, 1983.

\bibitem{Gr}
Phillip Griffiths, editor.
\newblock {\em {Topics in Transcendental Algebraic Geometry}}, volume 106 of
  {\em Ann. Math Studies}.
\newblock Princeton University Press, 1984.

\bibitem{Ha}
Yoshiko Hayakawa.
\newblock {Degeneration of Calabi-Yau Manifold with Weil-Petersson Metric}.
\newblock Technical Report alg-geom/9507016, Okolahoma State University, July
  1995.

\bibitem{S}
Andrew Strominger.
\newblock {Special Geometry}.
\newblock {\em Comm. Math. Phy.}, 133:163--180, 1990.

\bibitem{T1}
Gang Tian.
\newblock {Smoothness of the Universal Deformation Space of Compact Calabi-Yau
  Manifolds and its Peterson-Weil Metric}.
\newblock In Shing-Tung Yau, editor, {\em Mathematical aspects of string
  theory}, volume~1, pages 629--646. World Scientific, 1987.

\bibitem{T2}
Gang Tian.
\newblock {Smoothing 3-folds with Trivial Canonical Bundle and Ordinary Double
  Points}.
\newblock In Shing-Tung Yau, editor, {\em Eassays in Mirror Symmetry}, pages
  458--479. International Press, 1992.

\bibitem{Wong}
Chin-Long Wang.
\newblock {On the Weil-Petersson Metrics and Degeneration of Calabi-Yau
  Manifolds}.
\newblock Technical report, Harvard University, Dec. 1996.

\bibitem{Y1}
Shing-Tung Yau.
\newblock {On the Ricci Curvature of a Compact K\"ahler Manifold and the
  Complex Monge-Ampere Equation, I,}.
\newblock {\em Comm. Pure. Appl. Math}, 31:339--411, 1978.

\end{thebibliography}

\end{document}